# Symmetry Classification and Exact Solutions of a Variable Coefficient Nonlinear Space-Time Fractional Burgers' Equation

Manoj Gaur* and K. Singh**

## Abstract

In this paper, we investigate the symmetry properties of a variable coefficient nonlinear space-time fractional Burgers' equation. Fractional Lie symmetries and corresponding infinitesimal generators are obtained. With the help of the infinitesimal generators some group invariant solutions are deduced. Further, some exact solutions of fractional Burgers' equation are generated by the invariant subspace method.

**MSC 2010: Primary 26A33: Secondary 33E12, 34A08, 34K37, 35R11, 60G22**

**Keywords:** Fractional differential equations, Fractional integration, Fractional Lie group method.

# 1. Introduction

Partial differential equations (PDEs) have a wide range of applications in many fields, such as physics, engineering and chemistry which are fundamental for the mathematical formulation of continuum models [5]. In the $18^{th}$ century Euler, D'Alembert, Lagrange and Laplace started the study of PDEs as the key mode of analytical study of models in the physical science. After the work of Riemann during 19th century, PDEs turn out to be an essential tool in other branches of mathematics [2]. The Burgers' equation is a one dimensional non-linear partial differential equation, which is a simple form of the one dimensional Navier-Stokes equation. It was presented for the first time in a paper in 1940s from Burger. Later the Burgers' equation was studied by Cole [4] who gave a theoretical solution, based on Fourier series analysis, using the appropriate initial and boundary conditions. Burgers' equation has a large variety of applications in the modelling of water in unsaturated soil, dynamics of soil water, statistics of flow problems, mixing and turbulent diffusion, cosmology and seismology [1,4] and is considered as one of the most important PDE in the theory of nonlinear conservation laws. Burgers' equation is known as a model equation, and that is why it is important. This equation has different types and each of them has special applications.

Recently, fractional differential equations have found extensive applications in many fields. Many important phenomena in viscoelasticity, electromagnetics, material science, acoustics and electrochemistry are elegantly described with the help of fractional order differential equations. The fractional calculus has applications in many diverse areas, such as mathematical physics, viscoelasticity, transmission theory, reheology of soils, electric conductance of biological systems, modelling of neurons, diffusion processes, damping laws, and growth of intergranular grooves on metal surfaces [1,15,16]. It has been revealed that the non-conservative forces can be described by fractional differential equations. Therefore, as most of the processes in the real physical world are non-conservative, the fractional calculus can be used to describe them. Fractional integrals and derivatives also appear in the theory of control of dynamical systems, when the controlled system and/or the controller is described by a fractional differential equation. In the last few decades, the subject of the fractional calculus has caught the consideration of many researchers, who contributed to its development. Recently, some analytical and numerical methods [3,8] have been introduced to solve a fractional order differential equation. However, all the methods have an insufficient development as they allow one to find solutions only in case of linear equations and for isolated examples of nonlinear equations [1,15-19]. It is very well known that the Lie group method is the most effective technique in the field of applied mathematics to find exact solutions of ordinary and partial differential equations [10]. However, this approach is not yet applied much to investigate symmetry properties of fractional differential equations (FDEs). To the best of our knowledge there are only a few papers (for examples, [7,8,21] ) in which Lie symmetries and similarity solutions of some fractional differential equations have been discussed by some researchers. More recently, Jumarie [12] proposed the modified Riemann-Liouville derivative and Jumarie-Lagrange method [13], after which a generalized fractional characteristic method and a fractional Lie group method have been introduced by Wu [21,22] in order to solve a fractional order partial differential equation.

In this paper, we intend to apply Lie group method to solve a fractional order nonlinear Burgers' equation of the form

$$u_t^{(\alpha)} = f(t) u_x^{(2\beta)} + g(t)(u_x^{(\beta)})^2 \qquad x \in (0, \infty),\ t > 0,\ 0 < \alpha, \beta < 1 \quad (1.1)$$

which is a generalisation of the time-fractional nonlinear Burgers' equation examined by Wu [22]. This work is based on some basic elements of fractional calculus, with special emphasis on the modified Riemann-Liouville type.

The paper is organized as follows. In Sec. 2, we briefly describe some definitions and properties of fractional calculus. In Sec. 3, we obtain the symmetries for the Burgers' equation having six-dimensional Lie algebras. In Sec. 4, we analyse the reduced systems and find some of the invariant solutions of eqn. (1.1). Sec.5 contains application of invariant subspace method on fractional Burgers' equation (1.1). Finally, a conclusion is given in Sec. 6.

**2. Basic concepts in fractional calculus**

In this paper, the modified Riemann-Liouville derivative proposed by Jumarie [12] has been adopted. Some definitions are given which have been used throughout this work.

**2.1. Fractional Riemann-Liouville integral**

The fractional Riemann-Liouville integral of a continuous (but not necessarily differentiable) real valued function $f(x)$ with respect to $(dx)^\alpha$ is defined as [15,12]

$$_0I_x^\alpha f(x) = \frac{1}{\Gamma(\alpha)} \int_0^x (x-t)^{\alpha-1} f(t) dt = \frac{1}{\Gamma(1+\alpha)} \int_0^x f(t)(dt)^\alpha, \ 0 < \alpha \leq 1 \quad (2.1)$$

**2.2. Modified Riemann-Liouville derivative**

Through the fractional Riemann-Liouville integral, Jumarie [12] proposed the modified Riemann-Liouville derivative of $f(x)$ as

$$_0D_x^\alpha f(x) = \frac{1}{\Gamma(n-\alpha)} \frac{d^n}{dx^n} \times \int_0^x (x-t)^{n-\alpha-1} \big(f(t) - f(0)\big) dt, \ n-1 < \alpha < n. \quad (2.2)$$

**2.3. Some useful formulae**

Herein, we list some properties of modified Riemann-Liouville derivative which have been used in this paper

(i) $\quad df(x) = \dfrac{D_x^\alpha f(x)(dx)^\alpha}{\Gamma(1+\alpha)}$

(ii) $\quad D_x^\alpha (uv) = (D_x^\alpha u)v + u(D_x^\alpha v)$

(iii) $\quad D_t^\alpha f\{x(t)\} = \dfrac{df}{dx} D_t^\alpha x(t), \ 0 < \alpha < 1$, given $\dfrac{df}{dx}$ exist.

(iv) $\quad D_x^\alpha x^\beta = \dfrac{\Gamma(1+\beta)}{\Gamma(1+\beta-\alpha)} x^{\beta-\alpha}, \ 0 < \beta < 1$, where, $x^\beta$ is $\alpha$-differentiable.

(v) $\quad \int (dx)^\beta = x^\beta$

(vi) $\quad \Gamma(1+\beta)dx = (dx)^\beta$

The above formulae and details of their scope of applications and limitations can be found in [12,14].

**2.4. Characteristic method for fractional order differential equations**

Applying the fractional chain rule proposed by Jumarie [12]

$$du = \frac{\partial^\beta u(x,t)}{\Gamma(1+\beta)\partial x^\beta}(dx)^\beta + \frac{\partial^\alpha u(x,t)}{\Gamma(1+\alpha)\partial t^\alpha}(dt)^\alpha, \ 0<\alpha,\beta<1, \quad (2.3)$$

Wu [21,22] extended the characteristic method of first order linear partial differential equation to a linear fractional differential equation of the form

$$a(x,t)\frac{\partial^\beta u(x,t)}{\partial x^\beta} + b(x,t)\frac{\partial^\alpha u(x,t)}{\partial t^\alpha} = c(x,t), \ 0<\alpha,\beta<1 \quad (2.4)$$

and introduced the fractional characteristic equation [21]

$$\frac{(dx)^\beta}{\Gamma(1+\beta)a(x,t)} = \frac{(dt)^\alpha}{\Gamma(1+\alpha)b(x,t)} = \frac{du}{c(x,t)} \quad (2.5)$$

**3. Classification of symmetries of nonlinear space-time fractional Burgers' equation**

Herein, we investigate the symmetries and reductions of nonlinear space-time fractional Burgers' equation (1.1)

Let us assume that eqn. (1.1) admits the Lie symmetries of the form

$$\frac{\tilde{x}^{\beta}}{\Gamma(1+\beta)} = \frac{x^{\beta}}{\Gamma(1+\beta)} + \varepsilon\xi(x,t,u) + o(\varepsilon^2) \qquad (3.1a)$$

$$\frac{\tilde{t}^{\beta}}{\Gamma(1+\alpha)} = \frac{t^{\beta}}{\Gamma(1+\alpha)} + \varepsilon\tau(x,t,u) + o(\varepsilon^2) \qquad (3.1b)$$

$$\tilde{u} = u + \varepsilon\eta(x,t,u) + o(\varepsilon^2) \qquad (3.1c)$$

where $\varepsilon$ is the group parameter and $\xi, \tau$ and $\eta$ are the infinitesimals of the transformations for the independent and dependent variables, respectively.

The associated Lie algebra of infinitesimal symmetries of eqn. (1.1) is then the fractional vector field of the form

$$V = \xi(x,t,u)\frac{\partial^{\beta}}{\partial x^{\beta}} + \tau(x,t,u)\frac{\partial^{\alpha}}{\partial t^{\alpha}} + \eta(x,t,u)\frac{\partial}{\partial u} \qquad (3.2)$$

and its fractional second order prolongation is given by

$$pr^{(2)}V = \xi(x,t,u)\frac{\partial^{\beta}}{\partial x^{\beta}} + \tau(x,t,u)\frac{\partial^{\alpha}}{\partial t^{\alpha}} + \eta(x,t,u)\frac{\partial}{\partial u} + \eta^{t}\frac{\partial}{\partial u_{t}^{(\alpha)}} + \eta^{x}\frac{\partial}{\partial u_{x}^{(\beta)}} +$$
$$\eta^{xt}\frac{\partial}{\partial (u_{t}^{(\alpha)})_{x}^{(\beta)}} + \eta^{tt}\frac{\partial}{\partial u_{t}^{(2\alpha)}} + \eta^{xx}\frac{\partial}{\partial u_{x}^{(2\beta)}} \qquad (3.3)$$

Now for the invariance of eqn. (1.1) under the eqns. (3.1 a-c), we must have

$$pr^{(2)}V\left([\Delta u]\right) = 0, \qquad (3.4)$$

where $[\Delta u] = u_{t}^{(\alpha)} - f(t)u_{x}^{(2\beta)} - g(t)(u_{x}^{(\beta)})^2,$ \qquad (3.5)

or equivalently, if

$$\left(\eta^t - f(t)\eta^{xx} - 2g(t)u_x^{(\beta)}\eta^x - \tau f_t^{(\alpha)}u_t^{(2\beta)} - \tau g_t^{(\alpha)}(u_x^{(\beta)})^2\right)\Big|_{[\Delta u]=0} = 0. \qquad (3.6)$$

The generalised fractional prolongation vector fields $\eta^x, \eta^t$ and $\eta^{xx}$ are given by

$$\eta^x = \eta_x^{(\beta)} + u_x^{(\beta)}\eta_u - (\xi_x^{(\beta)} + u_x^{(\beta)}\xi_u)u_x^{(\beta)} + (\tau_x^{(\beta)} + u_x^{(\beta)}\tau_u)u_t^{(\alpha)} \qquad (3.7a)$$

$$\eta^t = \eta_t^{(\alpha)} + u_t^{(\alpha)}\eta_u - (\xi_t^{(\alpha)} + u_t^{(\alpha)}\xi_u)u_x^{(\beta)} + (\tau_t^{(\alpha)} + u_t^{(\alpha)}\tau_u)u_t^{(\alpha)} \qquad (3.7b)$$

$$\eta^{xx} = \eta_x^{(2\beta)} + u_x^{(\beta)}(\eta_u)_x^{(\beta)} - (\xi_x^{(2\beta)} + u_x^{(\beta)}(\xi_u)_x^{(\beta)})u_x^{(\beta)}$$

$$+ (\tau_x^{(2\beta)} + u_x^{(\beta)}(\tau_u)_x^{(\beta)})u_t^{(\alpha)}$$

$$+ u_x^{(\beta)}\left[(\eta_x^{(\beta)})_u + u_x^{(\beta)}\eta_{uu} - ((\xi_x^{(\beta)})_u + u_x^{(\beta)}\xi_{uu})u_x^{(\beta)} - ((\tau_x^{(\beta)})_u + u_x^{(\beta)}\tau_{uu})u_t^{(\alpha)}\right]$$

$$+ u_x^{(2\beta)}\left[\eta_u - (\xi_x^{(\beta)} + u_x^{(\beta)}\xi_u) - u_x^{(\beta)}\xi_u - \tau_u u_t^{(\alpha)}\right]$$

$$- (\xi_x^{(\beta)} + u_x^{(\beta)}\xi_u)u_x^{(2\beta)} + (\tau_x^{(\beta)} + u_x^{(\beta)}\tau_u)(u_t^{(\alpha)})_x^{\beta} - (u_t^{(\alpha)})_x^{\beta}(\tau_x^{(\beta)} + u_x^{(\beta)}\tau_u)$$

$$(3.7c)$$

Using eqn. (3.7) in eqn. (3.6) and equating the coefficient of various derivative terms to zero, we get the set of simplified determining equations as follows

$$\tau_u = 0 \qquad (3.8a)$$

$$\tau_x^{(\beta)} = 0 \qquad (3.8b)$$

$$\xi_u = 0 \tag{3.8c}$$

$$2g(t)\xi_x^{(\beta)} - g(t)\tau_t^{(\alpha)} - \tau g_t^{(\alpha)} - g(t)\eta_u - f(t)\eta_{uu} = 0 \tag{3.8d}$$

$$2f(t)\xi_x^{(\beta)} - f(t)\tau_t^{(\alpha)} - \tau f_t^{(\alpha)} = 0 \tag{3.8e}$$

$$f(t)\xi_x^{(2\beta)} - \xi_t^{(\alpha)} - 2g(t)\eta_x^{(\beta)} - 2f(t)(\eta_x^{(\beta)})_u = 0 \tag{3.8f}$$

$$\eta_t^{(\alpha)} - f(t)\eta_x^{(2\beta)} = 0 \tag{3.8g}$$

On solving eqn. (3.8g) by using fractional Lie group method we obtain a particular solution as

$$\eta = a_1 \frac{x^\beta}{\Gamma(1+\beta)} + a_2 F_t^\alpha + a_2 \frac{x^{2\beta}}{\Gamma(1+2\beta)} + a_3$$

where $F_t^{2\alpha} = f(t)$.

Using this value of $\eta$ in eqns. (3.8d) and (3.8e), we get

$$2\{g(t) - f(t)\}\xi_x^{(2\beta)} = 0$$

which brings forth the following possibilities; either

(i) $\xi_x^{(2\beta)} = 0$ or (ii) $f(t) = g(t)$

Case I: In this case, from the determining equations, we get

$$\xi = -2a_1 G(t) - 2a_2 G(t)\frac{x^\beta}{\Gamma(1+\beta)} + a_4 \frac{x^\beta}{\Gamma(1+\beta)} + a_5 \tag{3.9a}$$

$$\tau = \frac{1}{F_t^{(2\alpha)}}\left[-4a_2 H(t) + 2a_4 F_t^{(\alpha)} + a_6\right] \tag{3.9b}$$

$$\eta = a_1 \frac{x^\beta}{\Gamma(1+\beta)} + a_2 F_t^\alpha + a_2 \frac{x^{2\beta}}{\Gamma(1+2\beta)} + a_3 \quad (3.9c)$$

where $G_t^{(\alpha)}(t) = g(t)$, $H_t^{(\alpha)}(t) = F_t^{(2\alpha)}G(t)$, and $a_1, a_2, ..., a_6$ are six arbitrary constant parameters. Using the above obtained values in (3.8d) we also get $g(t) = kf(t)$, where $k$ is an arbitrary constant. This covers the case $f(t) = g(t)$. Further for $f(t) = g(t) = 1$ and $\beta = 1$ the infinitesimals can be reduced to those reported in [22], by setting the coefficients $a_1 = -c_5$, $a_2 = -2c_6$, $a_3 = c_3$, $a_4 = c_4$, $a_5 = c_1$, $a_6 = c_2$.

Hence, the fractional point symmetry generators admitted by the eqn. (1.1) are given by

$$V_1 = -2G(t)\frac{\partial^\beta}{\partial x^\beta} + \frac{x^\beta}{\Gamma(1+\beta)}\frac{\partial}{\partial u}$$

$$V_2 = -2G(t)\frac{x^\beta}{\Gamma(1+\beta)}\frac{\partial^\beta}{\partial x^\beta} - \frac{4H(t)}{F_t^{(2\alpha)}}\frac{\partial^\alpha}{\partial t^\alpha} + \left(F_t^\alpha + \frac{x^{2\beta}}{\Gamma(1+2\beta)}\right)\frac{\partial}{\partial u}$$

$$V_3 = \frac{\partial}{\partial u}$$

$$V_4 = \frac{x^\beta}{\Gamma(1+\beta)}\frac{\partial^\beta}{\partial x^\beta} + \frac{2F_t^{(\alpha)}}{F_t^{(2\alpha)}}\frac{\partial^\alpha}{\partial t^\alpha}$$

$$V_5 = \frac{\partial^\beta}{\partial x^\beta}, \text{ and}$$

$$V_6 = \frac{1}{F_t^{(2\alpha)}}\frac{\partial^\alpha}{\partial t^\alpha}$$

These infinitesimal generators can be used to determine a six parameter fractional Lie group of point transformation acting on $(x,t,u)$-space. It can be verified easily that the set $\{V_1, V_2, V_3, V_4, V_5, V_6\}$ forms a six dimensional Lie algebra under the Lie bracket $[X,Y] = XY - YX$ and its commutator table is as given below:

|       | $V_1$   | $V_2$       | $V_3$ | $V_4$  | $V_5$  | $V_6$       |
|-------|---------|-------------|-------|--------|--------|-------------|
| $V_1$ | 0       | 0           | 0     | $V_1$  | $-V_3$ | $2V_5$      |
| $V_2$ | 0       | 0           | 0     | $2V_2$ | $2V_1$ | $4V_4-2V_3$ |
| $V_3$ | 0       | 0           | 0     | 0      | 0      | 0           |
| $V_4$ | $-V_1$  | $-2V_2$     | 0     | 0      | $V_1$  | $2V_6$      |
| $V_5$ | $V_3$   | $-2V_1$     | 0     | $-V_5$ | 0      | 0           |
| $V_6$ | $-2V_5$ | $2V_3-4V_4$ | 0     | $-2V_6$| 0      | 0           |

The group transformation generated by the infinitesimal generators $V_i$, $(i=1,2,...,6)$ is obtained by solving the system of ordinary differential equations

$$\frac{(d\tilde{x})^\beta}{\Gamma(1+\beta)d\varepsilon} = \xi_i(\tilde{x},\tilde{t},\tilde{u}) \tag{3.10a}$$

$$\frac{(d\tilde{t})^\alpha}{\Gamma(1+\alpha)d\varepsilon} = \tau_i(\tilde{x},\tilde{t},\tilde{u}) \tag{3.10b}$$

$$\frac{d\tilde{u}}{d\varepsilon} = \eta_i(\tilde{x},\tilde{t},\tilde{u}) \tag{3.10c}$$

With the initial conditions

$$\tilde{x}\,|_{\varepsilon=0} = x$$

$$\tilde{t}\,|_{\varepsilon=0} = t \qquad (3.11)$$

$$\tilde{u}\,|_{\varepsilon=0} = u$$

Exponentiating the infinitesimal symmetries of eqn. (1.1), we get the one parameter groups $g_i(\varepsilon)$ generated by $V_i$ for $i = 1, 2, \ldots, 6$

$$g_1 : \left( \frac{x^\beta}{\Gamma(1+\beta)}, \frac{t^\alpha}{\Gamma(1+\alpha)}, u \right) \to \left( \frac{x^\beta}{\Gamma(1+\beta)} - \frac{2\varepsilon G}{\Gamma(1+\alpha)}, \frac{t^\alpha}{\Gamma(1+\alpha)}, u + \frac{\varepsilon x^\beta}{\Gamma(1+\beta)} \right)$$

$$g_2 : \left( \frac{x^\beta}{\Gamma(1+\beta)}, \frac{t^\alpha}{\Gamma(1+\alpha)}, u \right)$$

$$\to \left( \frac{x^\beta}{\Gamma(1+\beta)} e^{-2G\varepsilon}, \frac{G}{(1+2\varepsilon G)}, u + \frac{\varepsilon F_t^\alpha}{(1+2\varepsilon G)} + \frac{\varepsilon x^{2\beta}}{\Gamma(1+2\beta)} e^{-\frac{4\varepsilon G}{(1+2\varepsilon G)}} \right)$$

$$g_3 : \left( \frac{x^\beta}{\Gamma(1+\beta)}, \frac{t^\alpha}{\Gamma(1+\alpha)}, u \right) \to \left( \frac{x^\beta}{\Gamma(1+\beta)}, \frac{t^\alpha}{\Gamma(1+\alpha)}, u + \varepsilon \right)$$

$$g_4 : \left( \frac{x^\beta}{\Gamma(1+\beta)}, \frac{t^\alpha}{\Gamma(1+\alpha)}, u \right) \to \left( \frac{e^\varepsilon x^\beta}{\Gamma(1+\beta)}, \frac{e^{2\varepsilon} t^\alpha}{\Gamma(1+\alpha)}, u \right)$$

$$g_5 : \left( \frac{x^\beta}{\Gamma(1+\beta)}, \frac{t^\alpha}{\Gamma(1+\alpha)}, u \right) \to \left( \frac{x^\beta}{\Gamma(1+\beta)} + \varepsilon, \frac{t^\alpha}{\Gamma(1+\alpha)}, u \right)$$

$$g_6 : \left( \frac{x^\beta}{\Gamma(1+\beta)}, \frac{t^\alpha}{\Gamma(1+\alpha)}, u \right) \to \left( \frac{x^\beta}{\Gamma(1+\beta)}, \frac{t^\alpha}{\Gamma(1+\alpha)} + \varepsilon, u \right)$$

Now, since $g_i$ is a symmetry, if

$$u = f\left(\frac{x^\beta}{\Gamma(1+\beta)}, \frac{t^\alpha}{\Gamma(1+\alpha)}\right) \text{ is a solution of eqn. (1.1) the}$$

following $u_i$ are also solutions of eqn. (1.1)

$$u_1 = f\left(\frac{x^\beta}{\Gamma(1+\beta)} - \frac{2\varepsilon G}{\Gamma(1+\alpha)}, \frac{t^\alpha}{\Gamma(1+\alpha)}, u + \frac{\varepsilon x^\beta}{\Gamma(1+\beta)}\right) - \frac{\varepsilon x^\beta}{\Gamma(1+\beta)}$$

$$u_2 = f\left(\begin{array}{c}\dfrac{x^\beta}{\Gamma(1+\beta)} e^{-2G\varepsilon}, \dfrac{G}{(1+2\varepsilon G)}, \\ u + \dfrac{\varepsilon F_t^\alpha}{(1+2\varepsilon G)} + \dfrac{\varepsilon x^{2\beta}}{\Gamma(1+2\beta)} e^{-\frac{4\varepsilon G}{(1+2\varepsilon G)}}\end{array}\right) - \frac{\varepsilon F_t^\alpha}{(1+2\varepsilon G)} - \frac{\varepsilon x^{2\beta}}{\Gamma(1+2\beta)} e^{-\frac{4\varepsilon G}{(1+2\varepsilon G)}}$$

$$u_3 = f\left(\frac{x^\beta}{\Gamma(1+\beta)}, \frac{t^\alpha}{\Gamma(1+\alpha)}\right) - \varepsilon$$

$$u_4 = f\left(\frac{e^\varepsilon x^\beta}{\Gamma(1+\beta)}, \frac{e^{2\varepsilon} t^\alpha}{\Gamma(1+\alpha)}\right)$$

$$u_5 = f\left(\frac{x^\beta}{\Gamma(1+\beta)} + \varepsilon, \frac{t^\alpha}{\Gamma(1+\alpha)}\right)$$

$$u_6 = f\left(\frac{x^\beta}{\Gamma(1+\beta)}, \frac{t^\alpha}{\Gamma(1+\alpha)} + \varepsilon\right) \text{ and}$$

**4. Some exact solutions of the space-time fractional Burgers' equation**

In this section, we investigate some exact solutions of eqn. (1.1) corresponding to following infinitesimal generators.

(i) $V_1$

(ii) $V_4$

(iii) $V = nV_5 + mV_3$

(iv) $V = rV_5 + V_6$

(v) $V = sV_3 + V_6$

where $r$, $s$, $m$ and $n$ are arbitrary parameters.

**Theorem 4.1** Under the group of transformations $T(x,t) = \dfrac{t^\alpha}{\Gamma(1+\alpha)}$ and $\phi(T) = \dfrac{x^{2\beta}}{\Gamma(1+2\beta)} + 2G(t)u$ the Burgers' equation (1.1) reduces to a linear differential equation of first order $\phi'(T) - H_1(T)\phi(T) = -H_2(T)$, where $H_1(T) = \dfrac{G_t^{(\alpha)}(t)}{G(t)}$ and $H_2(T) = F_t^{(2\alpha)}$. Which admits a solution given by

$$u(x,t) = \dfrac{1}{2G(t)}\left[e^{\int H_1 dT}\left\{k_1 - \int H_2 e^{-\int H_1 dT} dT\right\} - \dfrac{x^{2\beta}}{\Gamma(1+2\beta)}\right],$$ where $k_1$ is an arbitrary constant.

**Proof:** Consider the infinitesimal generator $V_1$, given by

$$V_1 = -2G(t)\dfrac{\partial^\beta}{\partial x^\beta} + \dfrac{x^\beta}{\Gamma(1+\beta)}\dfrac{\partial}{\partial u} \qquad (4.1)$$

We find the resulting invariant solution by reducing eqn. (1.1) to a linear ordinary differential equation (4.5) using differential invariants. The fractional characteristic equation for $V_1$ is

$$\frac{\frac{dx^\beta}{\Gamma(1+\beta)}}{-2G(t)} = \frac{\frac{dt^\alpha}{\Gamma(1+\alpha)}}{0} = \frac{du}{\frac{x^\beta}{\Gamma(1+\beta)}} \qquad (4.2)$$

On solving the above fractional characteristic equations we obtain two functionally independent invariants as

$$T = \frac{t^\alpha}{\Gamma(1+\alpha)}, \text{ and } v = \frac{x^{2\beta}}{\Gamma(1+2\beta)} + 2G(t)u \qquad (4.3)$$

Now the solution of the fractional characteristic equation will be of the form $v = \phi(T)$, therefore,

$$u = \frac{1}{2G(t)}\left[\phi\left(\frac{t^\alpha}{\Gamma(1+\alpha)}\right) - \frac{x^{2\beta}}{\Gamma(1+2\beta)}\right] \qquad (4.4)$$

Substituting this value of u in eqn. (1.1), we get the reduced linear ode as

$$\frac{d\phi}{dT} - H_1(T)\phi = -H_2(T), \qquad (4.5)$$

where $T = \frac{t^\alpha}{\Gamma(1+\alpha)}, \quad H_1(T) = \frac{G_t^{(\alpha)}(t)}{G(t)}$ and $H_2(T) = F_t^{(2\alpha)}$

On solving eqn. (4.5), we obtain

$$\phi(T) = e^{\int H_1 dT} \left[ k_1 - \int H_2 e^{-\int H_1 dT} dT \right],$$

where $k_1$ is an arbitrary constant. This gives

$$u(x,t) = \frac{1}{2G(t)} \left[ e^{\int H_1 dT} \left\{ k_1 - \int H_2 e^{-\int H_1 dT} dT \right\} - \frac{x^{2\beta}}{\Gamma(1+2\beta)} \right] \quad (4.6)$$

**Theorem 4.2** The similarity transformations $u(x,t) = \psi(X)$ along with the similarity variable $X(x,t) = \dfrac{x^{2\beta}}{\left(\Gamma(1+\beta)\right)^2 F_t^{(\alpha)}}$ reduces the fractional Burgers' equation (1.1) to a nonlinear ordinary differential equation

$$\psi''(X) + \frac{g(t)}{f(t)}(\psi'(X))^2 + \frac{1}{4}\psi'(X) + \frac{1}{2X}\psi'(X) = 0 \text{ which leads to the}$$

solution $u(x,t) = \dfrac{1}{k} \log \left\{ k_2 + 2\pi K erf\left( \dfrac{x^{\beta}}{2\sqrt{F_t^{\alpha}}\Gamma(1+\beta)} \right) + \right\} + k_3$, where $k_2$

and $k_3$ are arbitrary constants.

**Proof:** Let us consider the infinitesimal generator

$$V_4 = \frac{x^{\beta}}{\Gamma(1+\beta)} \frac{\partial^{\beta}}{\partial x^{\beta}} + \frac{2F_t^{(\alpha)}}{F_t^{(2\alpha)}} \frac{\partial^{\alpha}}{\partial t^{\alpha}}. \quad (4.7)$$

Here, the fractional characteristic equations give the invariants

$$X(x,t) = \frac{x^{2\beta}}{\left(\Gamma(1+\beta)\right)^2 F_t^{(\alpha)}}, \text{ and } u(x,t) = \psi(X) \quad (4.8)$$

This yields

$$u(x,t) = \psi(X)$$

$$u = \psi\left(\frac{x^{2\beta}}{(\Gamma(1+\beta))^2 F_t^\alpha}\right) = \psi(X), \tag{4.9}$$

where $X = \dfrac{x^{2\beta}}{(\Gamma(1+\beta))^2 F_t^\alpha}$

Using the value of $u$ in eqn. (1.1), it becomes a nonlinear ordinary differential equation of second order

$$\psi''(X) + \frac{g(t)}{f(t)}(\psi'(X))^2 + \frac{1}{4}\psi'(X) + \frac{1}{2X}\psi'(X) = 0 \tag{4.10}$$

which leads to the solution

$$u(x,t) = \frac{1}{k}\log\left\{k_2 + 2\pi K erf\left(\frac{x^\beta}{2\sqrt{F_t^\alpha}\,\Gamma(1+\beta)}\right) + \right\} + k_3 \tag{4.11}$$

**Theorem 4.3** Under the transformations $\varsigma(x,t) = \dfrac{t^\alpha}{\Gamma(1+\alpha)}$ and $\varphi(\varsigma) = \dfrac{x^\beta}{\Gamma(1+\beta)} - \dfrac{n}{m}u$ the fractional Burgers' equation (1.1) reduces to an ordinary fractional differential equation $\varphi_t^{(\alpha)} + \dfrac{m}{n}G_t^{(\alpha)} = 0$ which has the general solution as $u(x,t) = \dfrac{m}{n}\left[\dfrac{x^\beta}{\Gamma(1+\beta)} + \dfrac{m}{n}G(t) - k_4\right]$, with $k_4$ as an arbitrary constant.

**Proof:** In this case, we study the infinitesimal generator

$$V = nV_5 + mV_3 = n\frac{\partial^\beta}{\partial x^\beta} + m\frac{\partial}{\partial u} \qquad (4.12)$$

The following invariants can be derived easily

$$\varsigma(x,t) = \frac{t^\alpha}{\Gamma(1+\alpha)}, v = \frac{x^\beta}{\Gamma(1+\beta)} - \frac{n}{m}u \qquad (4.13)$$

and the reduced form of eqn. (1.1) is

$$\varphi_t^{(\alpha)} + \frac{m}{n}G_t^{(\alpha)} = 0, \qquad (4.14)$$

This easily yields the solution

$$u(x,t) = \frac{m}{n}\left[\frac{x^\beta}{\Gamma(1+\beta)} + \frac{m}{n}G(t) - k_4\right] \qquad (4.15)$$

**Theorem 4.4** Under the transformations $\zeta(x,t) = \frac{1}{r}\frac{x^\beta}{\Gamma(1+\beta)} - \frac{F_t^{(\alpha)}}{\Gamma(1+\alpha)}$ and $\omega(\zeta) = u(x,t)$ the fractional Burgers' equation (1.1) reduces to a nonlinear ordinary differential equation $\omega''(\zeta) + k\omega'^2(\zeta) + \frac{r^2}{\Gamma(1+\alpha)}\omega'(\zeta) = 0$ which admits the solution

$$u(x,t) = ku = \frac{1}{k}\log\left[\left\{e^{\frac{r^2}{\Gamma(1+\alpha)}\left\{\frac{1}{r}\frac{x^\beta}{\Gamma(1+\beta)} - \frac{F_t^{(\alpha)}}{\Gamma(1+\alpha)}\right\}} - ke^{\frac{r^2}{\Gamma(1+\alpha)}c_1}\right\}\right] -$$

$$\frac{r^2}{\Gamma(1+\alpha)}\left\{\frac{1}{r}\frac{x^\beta}{\Gamma(1+\beta)} - \frac{F_t^{(\alpha)}}{\Gamma(1+\alpha)}\right\} + c_2,$$

where $c_1$ and $c_2$ are arbitrary constants.

**Proof**: In this case we obtain an invariant solution of eqn. (1.1) by using the infinitesimal generator

$$V = rV_5 + V_6 = r\frac{\partial^\beta}{\partial x^\beta} + \frac{1}{F_t^{(2\alpha)}}\frac{\partial^\alpha}{\partial t^\alpha}. \tag{4.16}$$

We have the invariants

$$\zeta(x,t) = \frac{1}{r}\frac{x^\beta}{\Gamma(1+\beta)} - \frac{F_t^{(\alpha)}}{\Gamma(1+\alpha)}, \text{ and } \omega(\zeta) = u(x,t) \tag{4.17}$$

and the reduced form of eqn. (1.1) as

$$\omega''(\zeta) + k\omega'^2(\zeta) + \frac{r^2}{\Gamma(1+\alpha)}\omega'(\zeta) = 0, \tag{4.18}$$

From eqn. (4.18), we get the solution as follows

$$u = \frac{1}{k}\log\left[\left\{e^{\frac{r^2}{\Gamma(1+\alpha)}\left\{\frac{1}{r}\frac{x^\beta}{\Gamma(1+\beta)} - \frac{F_t^{(\alpha)}}{\Gamma(1+\alpha)}\right\}} - ke^{\frac{r^2}{\Gamma(1+\alpha)}c_1}\right\}\right] - \frac{r^2}{\Gamma(1+\alpha)}\left\{\frac{1}{r}\frac{x^\beta}{\Gamma(1+\beta)} - \frac{F_t^{(\alpha)}}{\Gamma(1+\alpha)}\right\} + c_2 \tag{4.19}$$

**Theorem 4.5** Under the group of transformations $\gamma(x,t) = \dfrac{x^\beta}{\Gamma(1+\beta)}$ and $\rho(\gamma) = -\dfrac{\Gamma(1+\alpha)}{r} u + F_t^{(\alpha)}$ the Burgers' equation (1.1) reduces to a non-linear ordinary differential equation

$$\rho''(\gamma) - \dfrac{kr}{\Gamma(1+\alpha)} (\rho'(\gamma))^2 + 1 = 0.$$ Which has the general

solution as $u(x,t) = \left[ F_t^{(\alpha)} + \dfrac{1}{k} \log\left[ \cosh\left\{ \sqrt{k}(c_3 + \dfrac{x^\beta}{\Gamma(1+\beta)}) \right\} \right] - c_4 \right]$, where

$c_3$ and $c_4$ are arbitrary constants.

**Proof:** From the generator

$$V = sV_3 + V_6 = s \dfrac{\partial}{\partial u} + \dfrac{1}{F_t^{(2\alpha)}} \dfrac{\partial^\alpha}{\partial t^\alpha} \qquad (4.20)$$

we find the invariants

$$\gamma(x,t) = \dfrac{x^\beta}{\Gamma(1+\beta)}, \text{ and } \rho(\gamma) = -\dfrac{\Gamma(1+\alpha)}{s} u + F_t^{(\alpha)} \qquad (4.21)$$

and the reduced ode works out to be

$$\rho''(\gamma) - \dfrac{ks}{\Gamma(1+\alpha)} (\rho'(\gamma))^2 + 1 = 0, \qquad (4.22)$$

For the case under consideration the solution of eqn. (1.1) is

$$u(x,t) = \left[ F_t^{(\alpha)} + \frac{1}{k}\log\left[\cosh\left\{\sqrt{k}(c_3 + \frac{x^\beta}{\Gamma(1+\beta)})\right\}\right] - c_4 \right] \quad (4.23)$$

## 5. Some exact solutions of fractional Burgers' equation by the invariant subspace method.

The invariant subspace method was introduced by Galaktionov [6] in order to discover exact solutions of nonlinear partial differential equations. The method was further applied by Gazizov and Kasatkin [9] to some nonlinear fractional order differential equations. Here we give a brief description of the method.

Consider the fractional evolution equation $\frac{\partial^\alpha u}{\partial t^\alpha} = F[u]$, where $u = u(x,t)$ and $F[u]$ is a nonlinear differential operator.

The $n$-dimensional linear space $W_n = \langle f_1(x),...,f_n(x) \rangle$ is called invariant under the operator $F[u]$, iff $F[u] \in W_n$ for any $u \in W_n$. Which means there exist $n$ functions $\phi_1,...,\phi_n$ such that

$$F[C_1 f_1(x) + ... + C_n f_n(x)] = \phi_1(C_1,...,C_n) f_1(x) + ... + \phi_n(C_1,...,C_n) f_n(x),$$

where $C_1, C_2,..., C_n$ are arbitrary constants. The exact solution of fractional evolution equation can be obtained as $u(x,t) = \sum_{i=1}^{n} a_i(t) f_i(x)$.

For equation (1.1) $F[u] = f(t) u_x^{(2\beta)} + g(t)(u_x^{(\beta)})^2$. We have, the space

$$W_3 = \left\langle 1, \frac{x^\beta}{\Gamma(1+\beta)}, \frac{x^{2\beta}}{\Gamma(1+2\beta)} \right\rangle \text{ as invariant under } F[u], \text{ if and only if, for}$$

any $u \in W_n$, $F[u] \in W_3$.

Let $u(x,t) = a + b \dfrac{x^\beta}{\Gamma(1+\beta)} + c \dfrac{x^{2\beta}}{\Gamma(1+2\beta)}$ \hfill (5.1)

Then

$$F[u] = F\left[ a + b \frac{x^\beta}{\Gamma(1+\beta)} + c \frac{x^{2\beta}}{\Gamma(1+2\beta)} \right] = cf(t) + g(t)\left( b + c \frac{x^\beta}{\Gamma(1+\beta)} \right)^2.$$

Now $F[u] \in W_3$ if and only if,

$$cf(t) = a_1$$

$$g(t)b^2 = a_2$$

$$g(t)c^2 = a_3, \text{ and}$$

$$g(t)bc = a_4,$$

where $a_1, a_2, a_3,$ and $a_4$ are arbitrary constants.

Substituting the value of $u(x,t)$ from eqn. (5.1) into the eqn. (1.1), we get the following system of ordinary fractional differential equations

$$\frac{d^\alpha c}{dt^\alpha} = g(t)c^2$$

$$\frac{d^\alpha b}{dt^\alpha} = 2g(t)bc, \text{ and} \tag{5.2}$$

$$\frac{d^\alpha a}{dt^\alpha} = f(t)c + g(t)b^2$$

Eqns. (5.2) can be readily solved to yield

$$-\frac{1}{c} = \frac{1}{\Gamma(1+\alpha)} \int g(t)dt^\alpha + s_1$$

$$\log b = \frac{1}{\Gamma(1+\alpha)} \int 2cg(t)dt^\alpha + s_2, \text{ and}$$

$$a = \frac{1}{\Gamma(1+\alpha)} \left[ \int (cf(t) + b^2 g(t))dt^\alpha + s_3 \right],$$

where $s_1, s_2$ and $s_3$ are arbitrary constants.

## 6. Conclusion

The classical Lie group method is the most efficient technique for solving nonlinear partial differential equations. Here, in this paper, the fractional Lie group method has been effectively applied on a nonlinear space-time fractional Burgers' equation. For various infinitesimal generators the Burgers' equation has been reduced to some ordinary differential equations by using the method of differential invariants and some invariant solutions are obtained through the admitted symmetries and the invariant subspace method.

*Department of Mathematics, Jaypee University of Information Technology, Waknaghat, Solan 173234 (H.P.), INDIA

e-mail: *manojgaur22@gmail.com*

**Department of Mathematics, Jaypee University of Information Technology, Waknaghat, Solan 173234 (H.P.), INDIA

e-mail: *karan_jeet@yahoo.com*